\documentclass[12pt,twoside]{article}
\usepackage{amsmath}
\usepackage{epsfig}
\newtheorem{theor}{Theorem}
\newtheorem{lemma}[theor]{Lemma}

\newtheorem{examp}[theor]{Example}
\def\R{{\mathord{\rm I\mkern-3.6mu R}}}
\def\N{{\mathord{\rm I\mkern-3.6mu N}}}
\newcommand{\me}{1/2}

\newcommand{\vphi}{\varphi}
\newcommand{\lbd}{\lambda}
\newcommand{\D}{\protect\displaystyle}

\newcommand{\ipl}{\langle}
\newcommand{\ipr}{\rangle}

\begin{document}

\title{\centering  \large \bf ON INVERSE PROBLEMS MODELED BY PDE'S}
\author{\large \bf A.\,Leit\~ao%
\thanks{Research partially supported by a research grant from CNPq--GMD, 
\hfill \mbox{\hskip3cm}
{\it AMS Subject Classifiation}: Primary 65J20, Secondary 65L10, 65P05. 
\hfill \mbox{\hskip3cm}
{\it Key words and frases}:
Inverse problems, Reconstruction of boundary data, PDE's. }}
\date{}
\maketitle

\vspace{-5.6cm}
\noindent {\texttt{Proc.\,of\;VI\:Workshop\:on\:PDE's\;/\;%
July\:19\,-\,23,\:1999\;/\;Rio\:de\:Janeiro,\:Brazil}
\vspace{4cm}

\thispagestyle{empty}
\pagestyle{myheadings}
\markboth{A. LEIT\~AO}{ON INVERSE PROBLEMS MODELED BY PDE'S}

\begin{abstract}
We investigate the iterative methods proposed by Maz'ya and Kozlov 
(see [3], [4]) for solving ill-posed reconstruction problems modeled by PDE's. 
We consider linear time dependent problems of elliptic, hyperbolic and 
parabolic types. Each iteration of the analyzed methods consists on the 
solution of a well posed boundary (or initial) value problem. The iterations 
are described as powers of affine operators, as in [4]. We give alternative 
convergence proofs for the algorithms, using spectral theory and some 
functional analytical results (see [5], [6]).
\end{abstract} 

\centerline{\small\bf Resumo}
\begin{small}

        Investigamos neste artigo os m\'etodos propostos em [3] e [4] para 
resolver problemas mal postos de reconstru\c c\~ao. Consideramos problemas 
modelados por equa\c c\~oes elipticas, hiperb\'olicas e parab\'olicas. Cada 
itera\c c\~ao dos m\'etodos propostos se consiste na solu\c c\~ao de um 
problema bem posto. Apresentamos demonstra\c c\~oes alternativas as 
originais, utilizando argumentos de an\'alise funcional e teoria spectral 
(veja [5], [6]).
\end{small}

\setlength{\baselineskip}{7mm}
\section*{\large \bf 1. Introduction}

\subsection*{\normalsize \bf 1.1. Main results}

        We present new convergence proofs for iterative algorithms in [KM2]
using a functional analytical approach, where each iteration is described
using powers of an affine operator $T$. The key of the proof is to choose a
correct topology for the Hilbert space where the iteration takes place, and
to prove that $T_l$, the linear component of $T$, is a {\em regular
asymptotic}, {\em non expansive} operator.

         Other properties of $T_l$ such as positiveness, self-adjointness 
and injectivity are also verified. The ill-posed problems are presented in 
Section~2. In Section~3 we describe the iterative methods for each problem. 
The results concerning the analysis of the methods are summarized in 
Section~4. Some numerical results are discussed in Section~5.

        The iterative procedures discussed in this paper were first presented
in [KM2] and also treated in [Bas]. The iterative procedure for elliptic
(stationary) Cauchy problems is discussed in [KM1], [Le1,2] and [JoNa]. The
iterative procedure concerning parabolic problems is also treated in [Va].

\subsection*{\normalsize \bf 1.2. Preliminaries}

        Let $H$ be a separable Hilbert space endowed with inner product 
$\ipl \cdot,\cdot \ipr$ and norm $\| \cdot \|$. The operator $T: H \to H$ 
is said to be {\em regular asymptotic} in $x \in H$ if $\| T^{k+1}(x) - 
T^k(x) \| \to 0$ as $k\to\infty$. If the above property holds for every 
$x \in H$, we say that $T$ is regular asymptotic in $H$. The operator $T$ 
is called {\em non expansive} if $\|T\| \leq 1$. The next lemma is the key 
of the convergence proofs presented in this article.%
\footnote{A complete proof can be found in [Le1].}

\begin{lemma} \label{satz-nexp-asreg}

        Let $T:H \to H$ be a linear non expansive operator. With $\Pi$ we 
denote the orthogonal projector defined on $H$ onto $\ker(I-T)$. The following 
assertions are equivalent: \\[1ex]
a) $T$ is regular asymptotic in $H$; \\[1ex]
b) $\D\lim_{k\to\infty} T^k x = \Pi x$.
\end{lemma}

        Let $\Omega \subset \R^n$ be an open, bounded set with smooth 
boundary and $A$ be a positive, self-adjoint, unbounded operator (with 
discrete spectrum) densely defined on the Hilbert space $H := L^2(\Omega)$. 
Let $E_\lbd$, $\lbd \in \R$, denote the resolution of the identity 
associated to $A$. We construct the family of Hilbert spaces 
${\cal H}^s(\Omega)$, $s \geq 0$ as the domain of definition of the powers 
of $A$
\begin{equation} \label{def-Hs-raeume}
{\cal H}^s(\Omega) \ := \ \{ \vphi \in H \ | \ \| \vphi \|_s := \left(
\int_0^\infty (1+\lbd^2)^s d \ipl E_\lbd \vphi, \vphi \ipr \right)^{\me} 
  < \infty \} .
\end{equation}
The Hilbert spaces ${\cal H}^{-s}(\Omega)$ (with $s > 0$) are defined by 
duality:%
\footnote{Alternatively one can define ${\cal H}^{-s}(\Omega)$ as the 
completion of $H$ in the (-s)-norm defined in (\ref{def-Hs-raeume}).}
${\cal H}^{-s} := ({\cal H}^{s})'$. It follows directly from the definition 
that ${\cal H}^0(\Omega) = H$.

        An interesting case occurs when $A = (-\Delta)^{\me}$, where $\Delta$ 
is the Laplace--Beltrami operator on $\Omega$. In this particular case the 
identity ${\cal H}^s(\Omega) = H^{2s}_0(\Omega)$ holds, where 
$H^{s}_0(\Omega)$ is the Sobolev space of index $s$ according to Lions and 
Magenes (see [LiMa] pp. 54).

        Given $T>0$ we define the spaces $L^2(0,T; {\cal H}^s(\Omega))$ and 
$C(0,T; {\cal H}^s(\Omega))$ of functions $u: [0,T] \ni t \mapsto u(t) \in 
{\cal H}^s(\Omega)$. These are normed spaces if considered respectively 
with the norms
$$ \|u\|_{2;0,T;s} \ := \ 
   \left( \int_0^T \|u(t)\|_s^2 \, dt \right)^{\me} \ \ \ {\rm and} \ \ \ \
   \|u\|_{\infty;0,T;s} \ := \ \sup_{t\in[0,T]} \|u(t)\|_s . $$

\section*{\large \bf 2. The ill-posed problems}

\subsection*{\normalsize \bf 2.1. An elliptic problem}

        Given functions $(f,g) \in {\cal H}^{\me}(\Omega) \times 
{\cal H}^{-\me}(\Omega)$, find $u \in (V_e, \|\cdot\|_{V_e})$, where
\begin{eqnarray*}
V_e & := & \{ v \in L_2(0,T; {\cal H}^1(\Omega)) \ | \ 
   (\partial_t^2 - A^2)u = 0 \ {\rm in} \ (0,T) \times \Omega \} \\
\|u\|_{V_e} & := & \left(
      \int_0^T ( \|u(t)\|_1^2 + \|\partial_t u(t)\|_0^2 )\, dt \right)^{\me} ,
\end{eqnarray*}
that satisfies \\[2ex]
$ (P_e) \hfill \left\{  \begin{array}{l}
        (\partial_t^2 - A^2) u \, = \, 0\, ,\ {\rm in}\ (0,T) \times \Omega \\
        u(0)                   \, = \, f\, ,\ \ 
        \partial_t u(0)        \, = \, g\, .
       \end{array} \right. \hfill  $ \\[2ex]
Note that if $u \in V_e$, then $\partial_t u \in L_2(0,T; H)$ and adequate 
trace theorems (see [LiMa]) guarantee that $u(0), u(T) \in {\cal H}^{\me} 
(\Omega)$ and $\partial_t u(0), \partial_t u(T) \in {\cal H}^{-\me}(\Omega)$. 
The ill-posedness of problem $(P_e)$ can be easily verified from the explicit 
representation of it's solution:
\begin{equation}
  u(t,x) \ = \ \cosh(At) f(x) \ + \ \sinh(At) A^{-1} g(x) .
\end{equation}

\subsection*{\normalsize \bf 2.2. A hyperbolic problem}

        Given functions $f, g \in {\cal H}^1(\Omega)$ find $u \in (V_h, 
\|\cdot\|_{V_h})$, where
\begin{eqnarray*}
V_h & := & \{ v \in C(0,T; {\cal H}^1(\Omega)) \ | \ 
         \partial_t u \in C(0,T; H) \ {\rm and} \ 
         (\partial_t^2 + A^2)u = 0 \} \\
\|u\|_{V_h} & := & \sup_{t\in[0,T]} \left(
                        \|u(t)\|_1^2 + \|\partial_t u(t)\|_0^2 \right)^{\me} ,
\end{eqnarray*}
that satisfies \\[2ex]
$ (P_h) \hfill \left\{  \begin{array}{l}
        (\partial_t^2 + A^2) u \, = \, 0\, ,\ {\rm in}\ (0,T) \times \Omega \\
        u(0) \, = \, f \, ,\ \ 
        u(T) \, = \, g \, .
       \end{array} \right. \hfill  $ \\[2ex]
Note that if $u \in V_h$, then $u(0), u(T) \in {\cal H}^1(\Omega)$ and 
$\partial_t u(0), \partial_t u(T) \in H$. We assume further the numbers 
$k\pi/T$, $k=1,2,\dots$ are not eigenvalues of $A$.%
\footnote{If this condition is not satisfied, one can easily see that problem 
$(P_h)$ is not uniquely solvable.}
This hyperbolic (Dirichlet) boundary value problem is ill-posed if the 
distance from the set $M := \{ k\pi/T ; k \in \N \}$ to $\sigma(A)$ (the 
spectrum of $A$) is zero. This follows from the explicit representation of 
the solution of $(P_h)$
\begin{equation}
  u(t,x) \ = \ \sin(A(T-t)) \sin(AT)^{-1} f(x) \ + \ 
                 \sin(At) \sin(AT)^{-1} g(x) .
\end{equation}

\subsection*{\normalsize \bf 2.1. A parabolic problem}

        Given a function $f \in H = L_2(\Omega)$ find $u \in (V_p, 
\|\cdot\|_{V_p})$, where
\begin{eqnarray*}
V_p & := & \{ v \in L_2(0,T; {\cal H}^1(\Omega)) \ | \ 
         (\partial_t + A^2)u = 0 \ {\rm in}\ (0,T) \times \Omega \} \\
\|u\|_{V_p} & := & \left(
   \int_0^T ( \|u(t)\|_1^2 + \|\partial_t u(t)\|_{-1}^2 )\, dt \right)^{\me} ,
\end{eqnarray*}
that satisfies \\[2ex]
$ (P_p) \hfill \left\{  \begin{array}{l}
          (\partial_t + A^2) u \, = \, 0\, ,\ {\rm in}\ (0,T) \times \Omega \\
          u(T)               \, = \, f\, .
       \end{array} \right. \hfill  $ \\[2ex]
Note that if $u \in V_p$, then $u(0), u(T) \in H$. This corresponds to the 
well known problem of inverse heat transport, which is known to be severely 
ill-posed. The solution of $(P_p)$ has the explicit representation
\begin{equation} \label{loesung-parab-probl}
  u(t,x) \ = \ \exp(A^2(T-t)) f(x) .
\end{equation}

\section*{\large \bf 3. Description of the Methods}

\subsection*{\normalsize \bf 3.1. An iterative procedure for the elliptic 
problem}

        Consider problem $(P_e)$ with data $(f,g) \in {\cal H}^{\me}(\Omega) 
\times {\cal H}^{-\me}(\Omega)$. Given any initial guess $\vphi_0 \in 
{\cal H}^{-\me}(\Omega)$ for $\partial_t u(T)$ we try to improve it by 
solving the following mixed boundary value problems (BVP) of elliptic type
$$ \left\{  \begin{array}{l}
     (\partial_t^2 - A^2) v \, = \, 0\, ,\ {\rm in}\ (0,T) \times \Omega \\
     v(0) \, = \, f\, ,\ \ \partial_t v(T) \, = \, \vphi
   \end{array} \right. \hskip1cm
   \left\{  \begin{array}{l}
     (\partial_t^2 - A^2) w \, = \, 0\, ,\ {\rm in}\ (0,T) \times \Omega \\
     \partial_t w(0) \, = \, g\, ,\ \ w(T) \, = \, v(T)
   \end{array} \right. $$
and defining $\vphi_1 := \partial_t w(T)$. Each one of the mixed BVP's 
above has a solution in $V_e$ and consequently $\vphi_1 \in {\cal H}^{-\me} 
(\Omega)$. Repeating this procedure we can construct a sequence $\{\vphi_k\}$ 
in ${\cal H}^{-\me}(\Omega)$. Using the explicit representation of the 
solutions $v$ and $w$ of the above problems, one obtains
$$ \vphi_1(x) \ = \
   \tanh(AT)^2 \vphi(x) \, + \, \sinh(At) \cosh(AT)^{-2} \, A f(x)
   \, + \, \cosh(AT)^{-1} \, g(x) . $$
Defining the affine operator $T_e: {\cal H}^{-\me} \to {\cal H}^{-\me}$, 
$T_e(\vphi) := \tanh(AT)^2 \vphi + h_{f,g}$, where $h_{f,g} := \sinh(At) 
\cosh(AT)^{-2} A f + \cosh(AT)^{-1} g$, the iterative algorithm can be 
rewritten as
\begin{equation}
\vphi_k \ = \ T_e(\vphi_{k-1}) \ = \ T_e^k(\vphi_0) \ = \ \tanh(AT)^{2k} 
\vphi_0 + \sum_{j=0}^{k-1} \tanh(AT)^{2j} h_{f,g} .
\label{elliptic-iter-def} \end{equation}

\subsection*{\normalsize \bf 3.2. An iterative procedure for the hyperbolic 
problem}

        Let's now consider problem $(P_h)$ with data $f, g \in {\cal H}^1 
(\Omega)$. Given any initial guess $\vphi_0 \in H$ for $\partial_t u(0)$ we 
try to improve it by solving the following initial value problems (IVP) of 
hyperbolic type%
\footnote{The second problem is considered with reversed time.}
$$ \left\{  \begin{array}{l}
     (\partial_t^2 + A^2) v \, = \, 0\, ,\ {\rm in}\ (0,T) \times \Omega \\
     v(0) \, = \, f\, ,\ \ \partial_t v(0) \, = \, \vphi
   \end{array} \right. \hskip1cm
   \left\{  \begin{array}{l}
     (\partial_t^2 + A^2) w \, = \, 0\, ,\ {\rm in}\ (0,T) \times \Omega \\
     w(T) \, = \, g\, ,\ \ \partial_t w(T) \, = \, \partial_t v(T)
     
   \end{array} \right. $$
and defining $\vphi_1 := \partial_t w(0)$. Each one of the IVP's above has a 
solution in $V_h$ and consequently $\vphi_1 \in H$. Repeating this procedure 
we can construct a sequence $\{\vphi_k\}$ in $H$. Determining the solutions 
$v$ and $w$ of the above problems, one obtains
$$ \vphi_1(x) \, = \, \partial_t w(0,x) \, = \, \cos(AT)^2 \vphi(x) -
                   \cos(AT)\sin(AT)\, A f(x) + \sin(AT)\, g(x) . $$
Now defining the affine operator $T_h: H \to H$, $T_h(\vphi) := \cos(AT)^2 
\vphi \ + \ h_{f,g}$, where $h_{f,g} := -\cos(AT)\sin(AT) A f + \sin(AT) g$, 
the iteration can be written as
\begin{equation}
\vphi_k \, = \, T_h(\vphi_{k-1}) \, = \, T_h^k(\vphi_0) \, = \,
\cos(AT)^{2k} \vphi_0 + \sum_{j=0}^{k-1} \cos(AT)^{2j} h_{f,g}.
\label{hyperbolic-iter-def} \end{equation}

\subsection*{\normalsize \bf 3.3. An iterative procedure for the parabolic 
problem}

        We consider problem $(P_p)$ with data $f \in H$. Let $\vphi_0 \in H$ 
be an  initial guess for $u(0)$ and define $\bar{\lbd} := \inf\{ \lbd; \lbd 
\in \sigma(A) \}$. Now choose a positive parameter $\gamma$ such that 
$\gamma < 2 \exp(\bar{\lbd}^2 T)$. The method consists in solving the 
initial value problems of parabolic type
$$ \left\{  \begin{array}{l}
     (\partial_t + A^2) v_0 \, = \, 0\, ,\ {\rm in}\ (0,T) \times \Omega \\
     v_0(0) \, = \, \vphi_0
   \end{array} \right. \hskip0.6cm
   \left\{  \begin{array}{l}
     (\partial_t + A^2) v_k \, = \, 0\, ,\ {\rm in}\ (0,T) \times \Omega \\
     v_k(0) \, = \, v_{k-1}(0) - \gamma (v_{k-1}(T) - f)
   \end{array} \right. $$
for $k \geq 1$. The sequence $\{\vphi_k\}$ is now defined by $\vphi_k := 
v_k(0) \in H$. Determining the solutions $v_k$ of the above problems, we 
have $\vphi_{k+1}(x) = (I - \gamma \exp(-A^2 T)) \vphi_k(x) + \gamma f(x)$. 
Define the affine operator $T_p: H \to H$, $T_p(\vphi) := 
(I - \gamma \exp(-A^2 T))\, \vphi \ + \ h_f$, where $h_f := \gamma f$, we 
can write the iteration as
\begin{eqnarray}
\vphi_k & = & T_p(\vphi_{k-1}) \ = \ T_p^k(\vphi_0) \nonumber \\ 
        & = & (I - \gamma \exp(-A^2 T))^k \vphi_0 \, + \, \sum_{j=0}^{k-1}
              (I - \gamma \exp(-A^2 T))^j h_f .
\label{parabolic-iter-def} \end{eqnarray}

\section*{\large \bf 4. Analysis of the Methods}

\subsection*{\normalsize \bf 4.1. The elliptic case}

        We start presenting a result, which is a generalization of the 
Cauchy--Kowalews\-ki theorem. A complete proof can be found in [Le1].

\begin{lemma}  Given $(f,g) \in {\cal H}^{\me} \times {\cal H}^{-\me}$, the 
problem $(P_e)$ has at most one solution in $V_e$. \hskip.5cm \qquad
\label{lemma-cauchy-kowal} \end{lemma}

        In the next theorem we verify some properties of the operator 
$T_{l,e}$, that will be needed for the convergence proof of the algorithm.

\begin{theor}  The linear operator $T_{l,e}$ is positive, self-adjoint, 
injective, regular asymptotic, non-expansive and $1 \not\in \sigma_p 
(T_{l,e})$. \\
\begin{proof}  The injectivity follows from 
Lemma~\ref{lemma-cauchy-kowal}. The properties: positiveness, 
self-adjointness and $1 \not \in \sigma_p(T_{l,e})$ follow from the 
assumptions on $A$. The last two properties follow from the inequality
\begin{equation} \label{gl-mazya-bedingung}
  \| (I-T_{l,e})x \|^2 \ \leq \ ( \| x \|^2 - \| T_{l,e} x \|^2 ) ,\
  \forall x \in {\cal H}^{-\me} . \qquad
\end{equation}
\end{proof}
\label{satz-Tle-eigenschaft} \end{theor}

\begin{theor}  If problem $(P_e)$ is consistent for the data $(f,g)$, then 
the sequence $\vphi_k$ converges to $\partial_t u(T)$ in the norm of 
${\cal H}^{-\me}(\Omega)$. \\
\begin{proof}  Follows from Theorem~\ref{satz-Tle-eigenschaft} and 
Lemma~\ref{satz-nexp-asreg}. \qquad
\end{proof}
\label{satz-converg-ellip-verfahren} \end{theor}

        The converse of Theorem~\ref{satz-converg-ellip-verfahren} is also 
true. This is discussed in

\begin{theor}  If the sequence $\vphi_k$ converges, say to $\bar{\vphi}$, 
then problem $(P_e)$ is consistent for the data $(f,g)$ and it's solution 
$u \in V_e$ satisfies $\partial_t u(T) = \bar{\vphi}$. \\
\begin{proof}  Follows from Lemma~\ref{lemma-cauchy-kowal} and the 
definition of $T_e$. \qquad
\end{proof}
\label{conversesatz-ellip} \end{theor}

\subsection*{\normalsize \bf 4.2. The hyperbolic case}

\begin{theor}  The linear operator $T_{l,h}: H \to H$ is positive, 
self-adjoint, injective, non-expansive, regular asymptotic and 1 is not an 
eigenvalue of $T_{l,h}$. \\
\begin{proof}  Analog to the proof of Theorem~\ref{satz-Tle-eigenschaft}.
\qquad
\end{proof}
\label{satz-Tlh-eigenschaft} \end{theor}

\begin{theor}  If problem $(P_h)$ is consistent for the data $(f,g)$, then 
the sequence $\vphi_k$ converges to $\partial_t u(0)$ in the norm of $H$. \\
\begin{proof}  Follows from Theorem~\ref{satz-Tlh-eigenschaft} and 
Lemma~\ref{satz-nexp-asreg}. \qquad
\end{proof}
\end{theor}

\begin{theor}  If the sequence $\vphi_k$ converges, say to $\bar{\vphi}$, 
then problem $(P_h)$ is consistent for the Cauchy data $(f,g)$ and it's 
solution $u \in V_h$ satisfies $\partial_t u(0) = \bar{\vphi}$. \\
\begin{proof}  Analog to the proof of Theorem~\ref{conversesatz-ellip}. \qquad
\end{proof}
\end{theor}

\subsection*{\normalsize \bf 4.3. The parabolic case}

\begin{lemma}  Given $f \in H$, the problem $(P_p)$ has exactly one 
solution in $V_p$. \\
\begin{proof}  This result is suggested by the general representation of the 
solution given in (\ref{loesung-parab-probl}). A complete proof can be found 
in [LiMa], Chapter~3. \qquad
\end{proof}
\end{lemma}

\begin{theor}  The linear operator $T_{l,p}: H \to H$ is self-adjoint, 
non-expansive, regular asymptotic and 1 is not an eigenvalue of $T_{l,p}$. 
Further, if it is possible to choose $\gamma <  2\exp(\tilde{\lbd}^2T)$, 
where $\tilde{\lbd} := ( \bar{\lbd}^2 - T^{-1} \ln 2 )^{\me}$, then 
$T_{l,p}$ is also injective. \\
\begin{proof}  Analog to the proof of Theorem~\ref{satz-Tle-eigenschaft}. 
The injectivity under the extra assumption on $\gamma$ follows from an 
inequality similar to (\ref{gl-mazya-bedingung}). \qquad
\end{proof}
\label{satz-Tlp-eigenschaft} \end{theor}

\begin{theor}  Given $f \in H$, let $u \in V_p$ be the uniquely determined 
solution of problem $(P_p)$. Then the sequence $\vphi_k$ converges to 
$u(0)$ in the norm of $H$. \\
\begin{proof}  Follows from Theorem~\ref{satz-Tlp-eigenschaft} and 
Lemma~\ref{satz-nexp-asreg}. \qquad
\end{proof}
\end{theor}

\section*{\large \bf 5. Numerical results}

%
\begin{examp}  Consider the problem of finding $u(0) \in L^2(\Omega)$, where 
$u$ solves
$$  \left\{ \begin{array}{r@{\ }c@{\ }l}
      a^2 \, \partial_t u \, - \, \Delta u & = & 0 \\ u(T) & = & f
    \end{array} \right. $$
In this example $\Omega = [0,1]\times[0,1]$, $a^2=2$ and the final time is 
$T=0.625$. We choose the parameter $\gamma=2$ for the iteration. In Figure~1 
one can see the problem data $f$ and the corresponding solution $u(0)$.

        In Figure~2 the error $| \vphi_k - u(0) |$ is is shown after $10$, 
$10^4$, $10^5$ and $10^6$ iterations. One should note that the reconstruction 
error is smaller at the part of the domain where $u(0)$ is smooth. In Table~1 
the evolution of the relative error in the $L^2$--norm of the iteration is 
shown. Note that the convergence speed decays exponentially as we iterate.

\centerline{ \epsfysize6cm \epsfxsize7cm \epsfbox{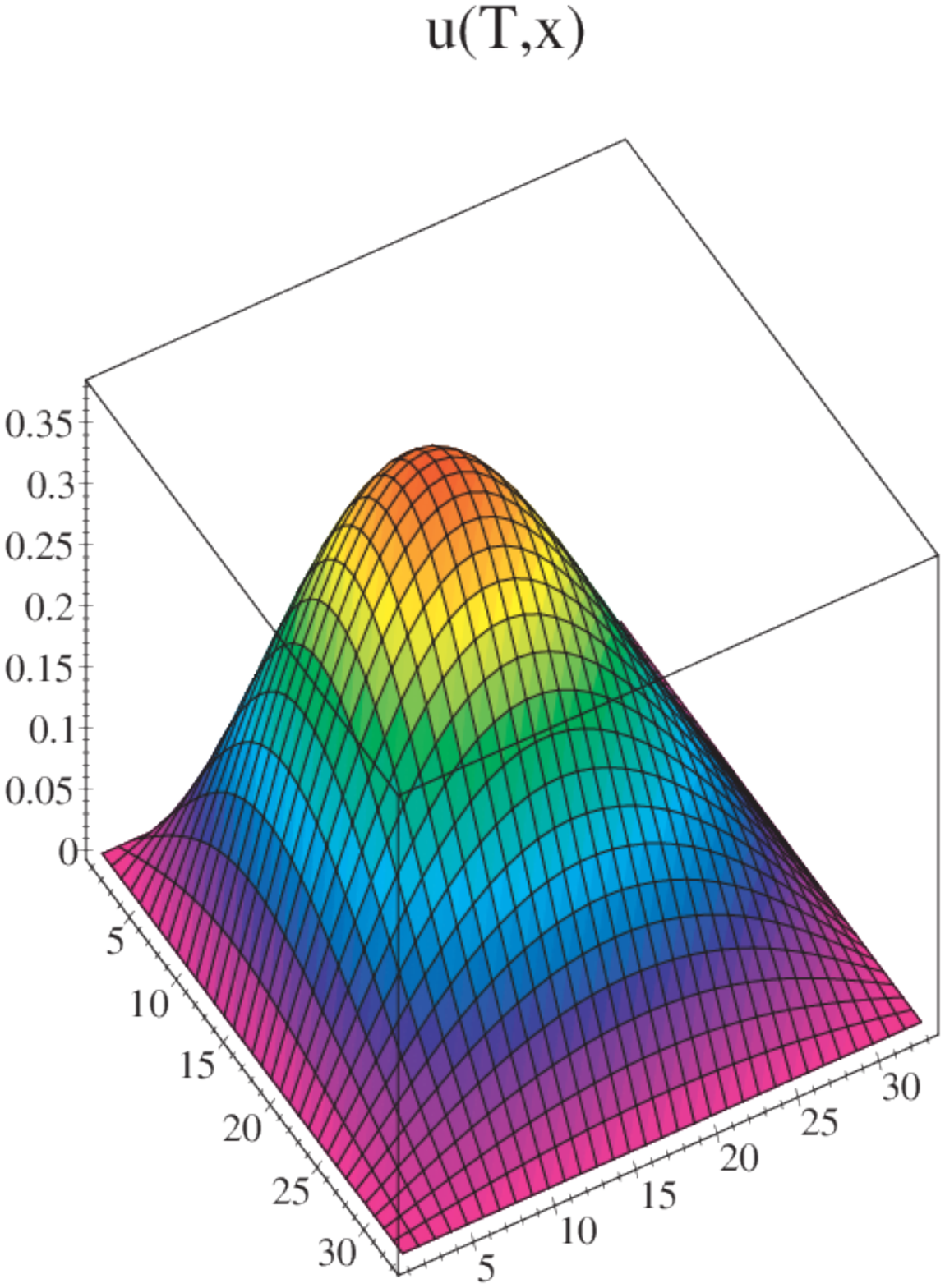}
             \epsfysize6cm \epsfxsize7cm \epsfbox{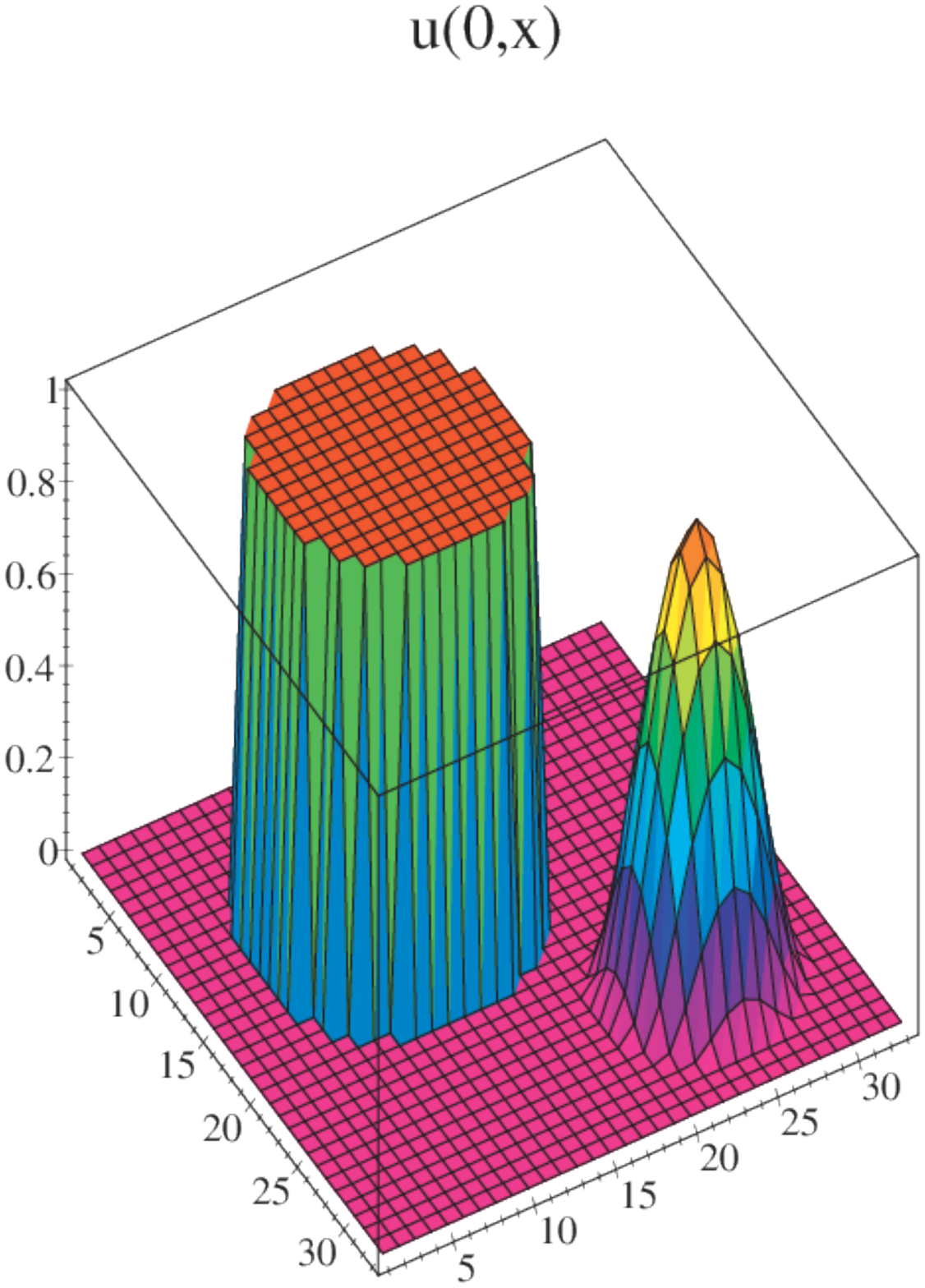} }

{\bf Figure 1:} Problem data $u(T)=f$ and corresponding $u(0)$.
\bigskip \bigskip

\centerline{ \epsfysize6cm \epsfxsize7cm \epsfbox{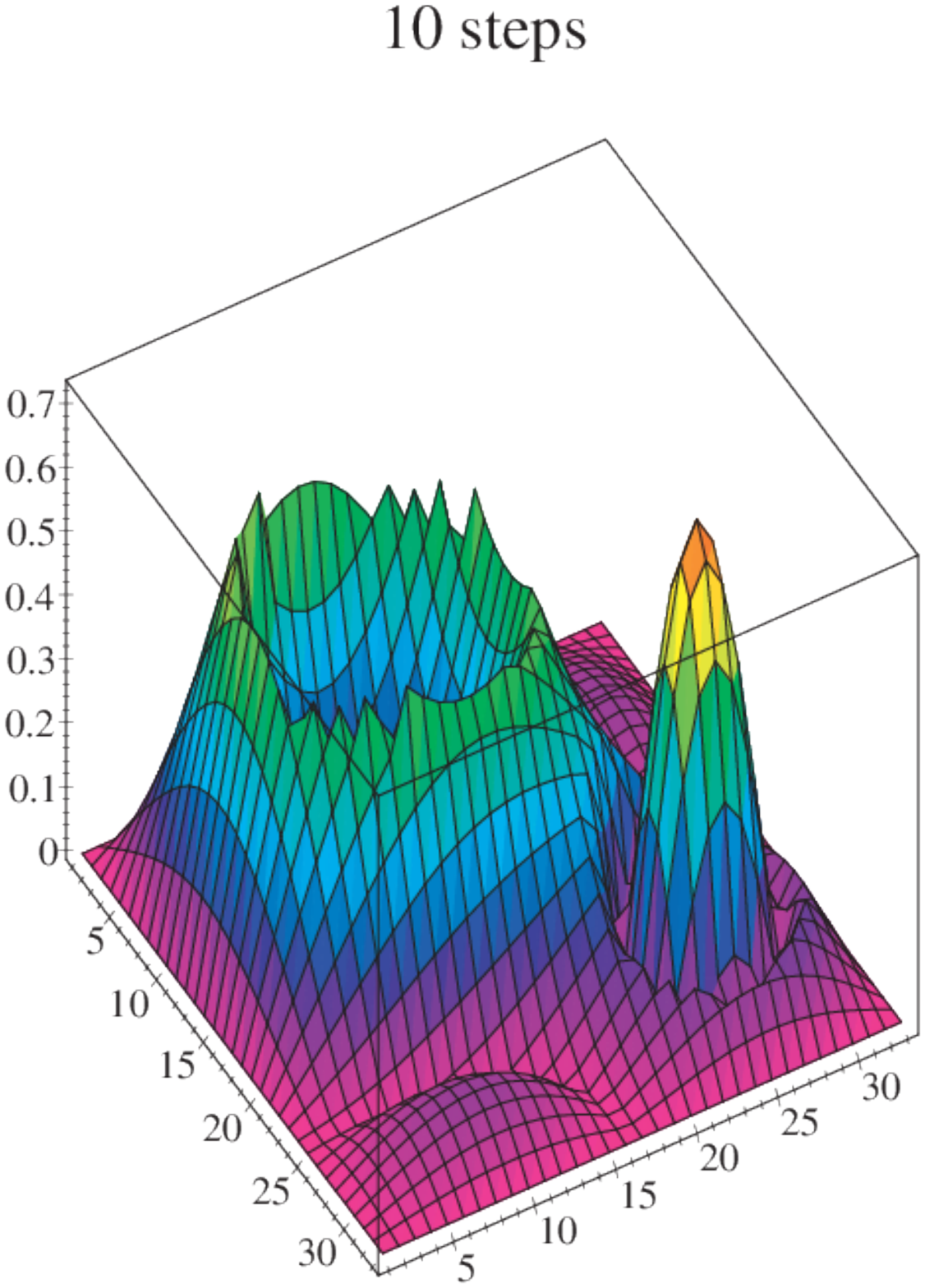}
             \epsfysize6cm \epsfxsize7cm \epsfbox{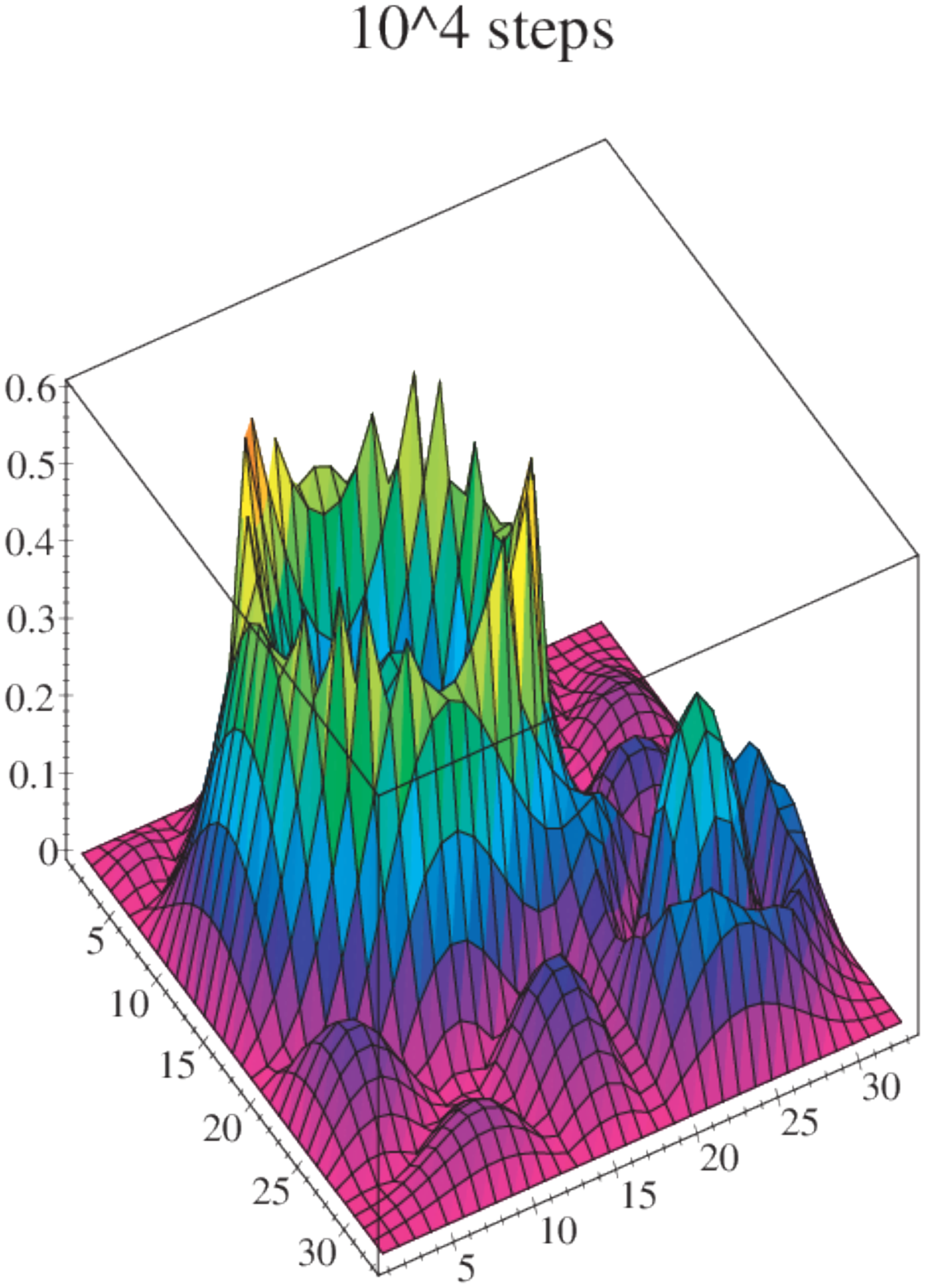} }
\vskip0.5cm

\centerline{ \epsfysize6cm \epsfxsize7cm \epsfbox{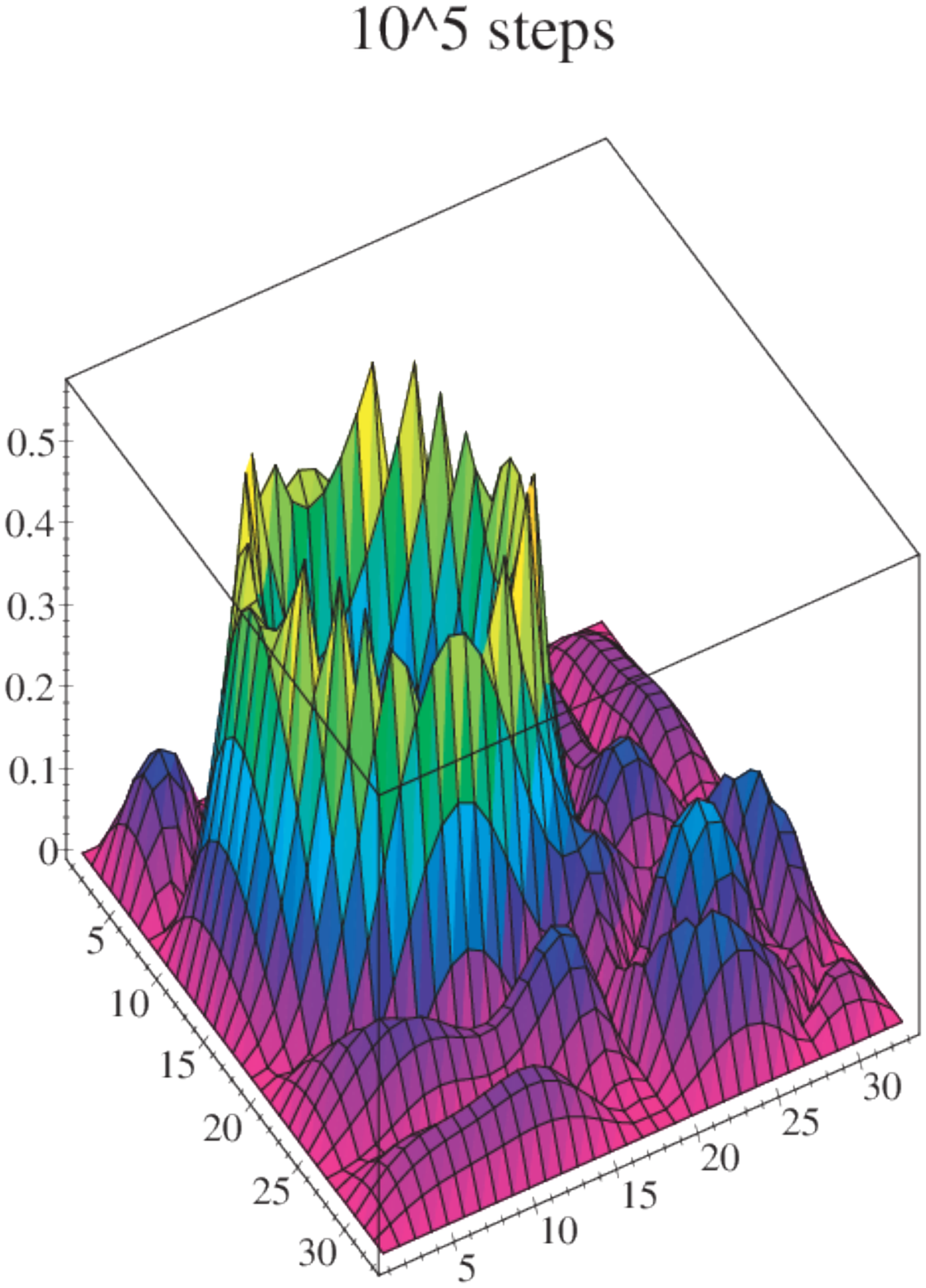}
             \epsfysize6cm \epsfxsize7cm \epsfbox{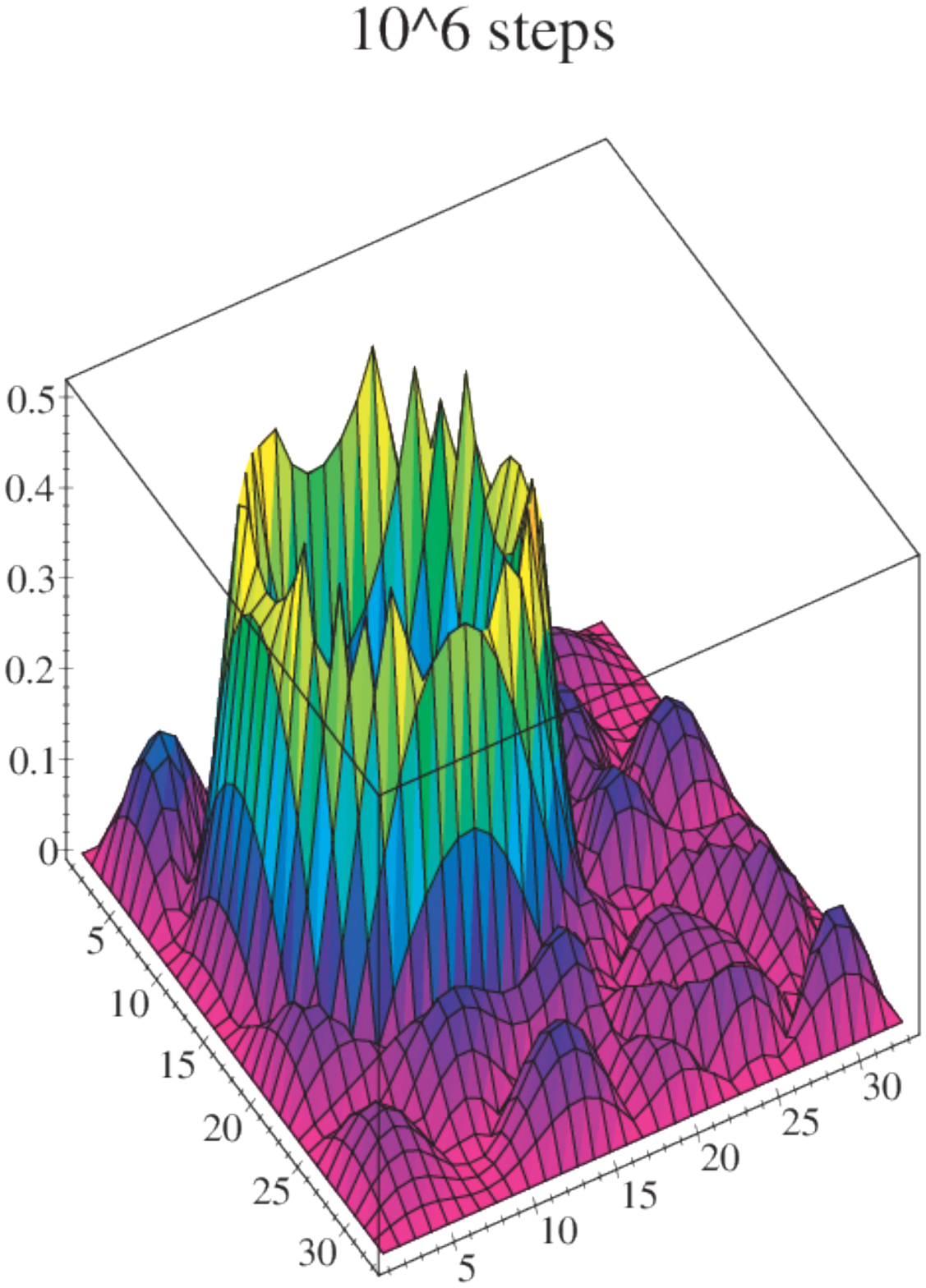} }

{\bf Figure 2:} Evolution of the error $| \vphi_k - u(0) |$.
\newpage

\begin{center} \begin{tabular}{@{}cccccc} \hline\hline
 $10$ steps & $10^3$ steps & $10^4$ steps & $10^5$ steps & $10^6$ steps \\
\hline
 49.8\%  & 42.2\% & 40.1\% & 36.2\% & 31.4\%  \\ \hline\hline
\end{tabular} \\[1ex]
{\bf Table 1:} Evolution of the relative error in the $L^2$--norm.
\end{center}

\end{examp}

%
\begin{examp}  Consider the problem of finding $w(0) \in L^2(\Omega)$, where 
$w$ solves:
$$  \left\{ \begin{array}{r@{\ }c@{\ }l}
       \partial_t w \, - \, \Delta w & = & 0 \\ w(T) & = & f
    \end{array} \right. $$
The set $\Omega$ is the same as in the previous example, the final time is
$T=0.625$ and we chose the parameter $\gamma=2$ for the iteration. We solve 
the direct Cauchy-problem for two different initial conditions, which are 
shown in Figures~3 and 4 respectively.
\bigskip

\centerline{ \epsfysize6cm \epsfxsize7cm \epsfbox{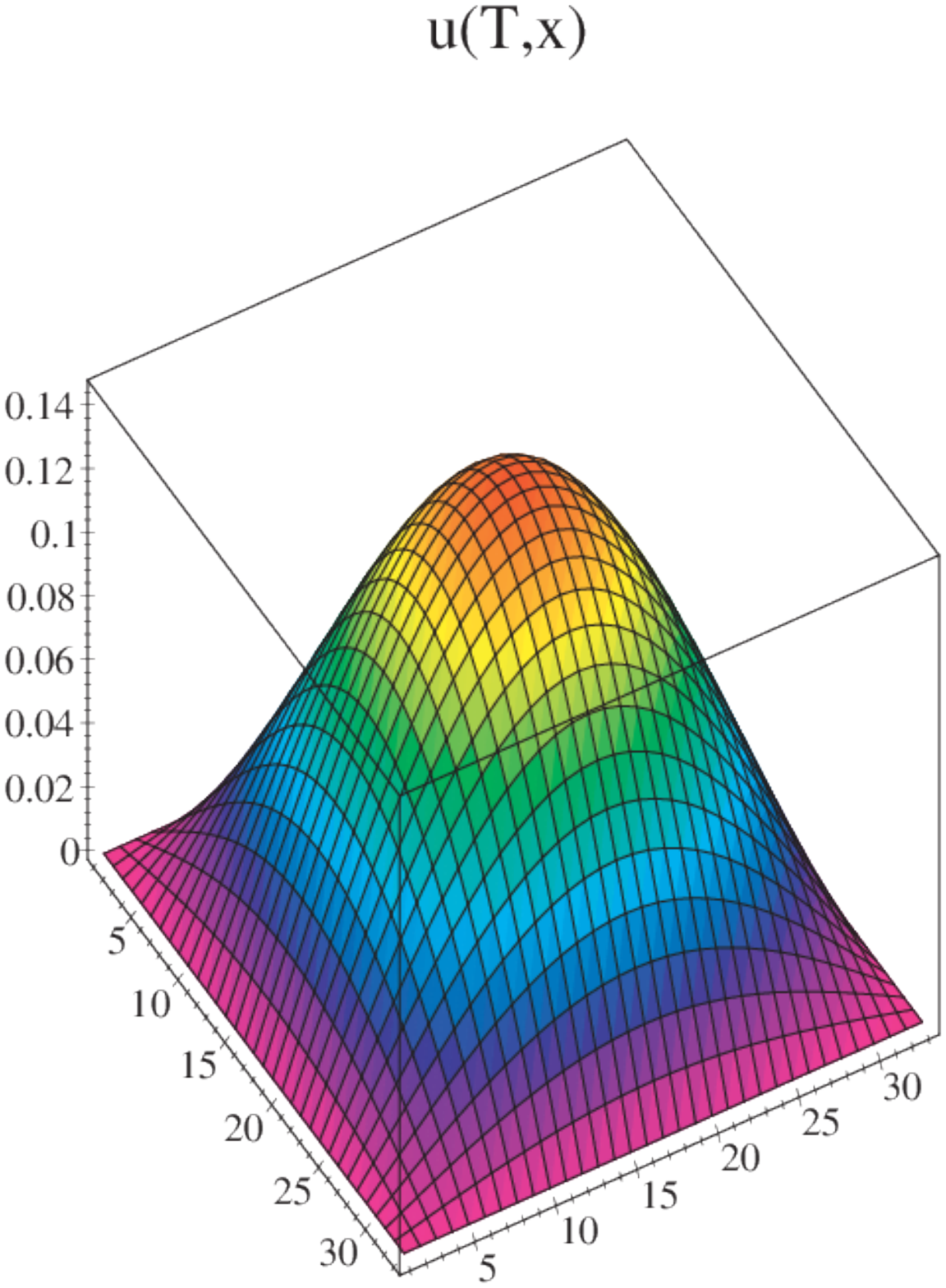}
             \epsfysize6cm \epsfxsize7cm \epsfbox{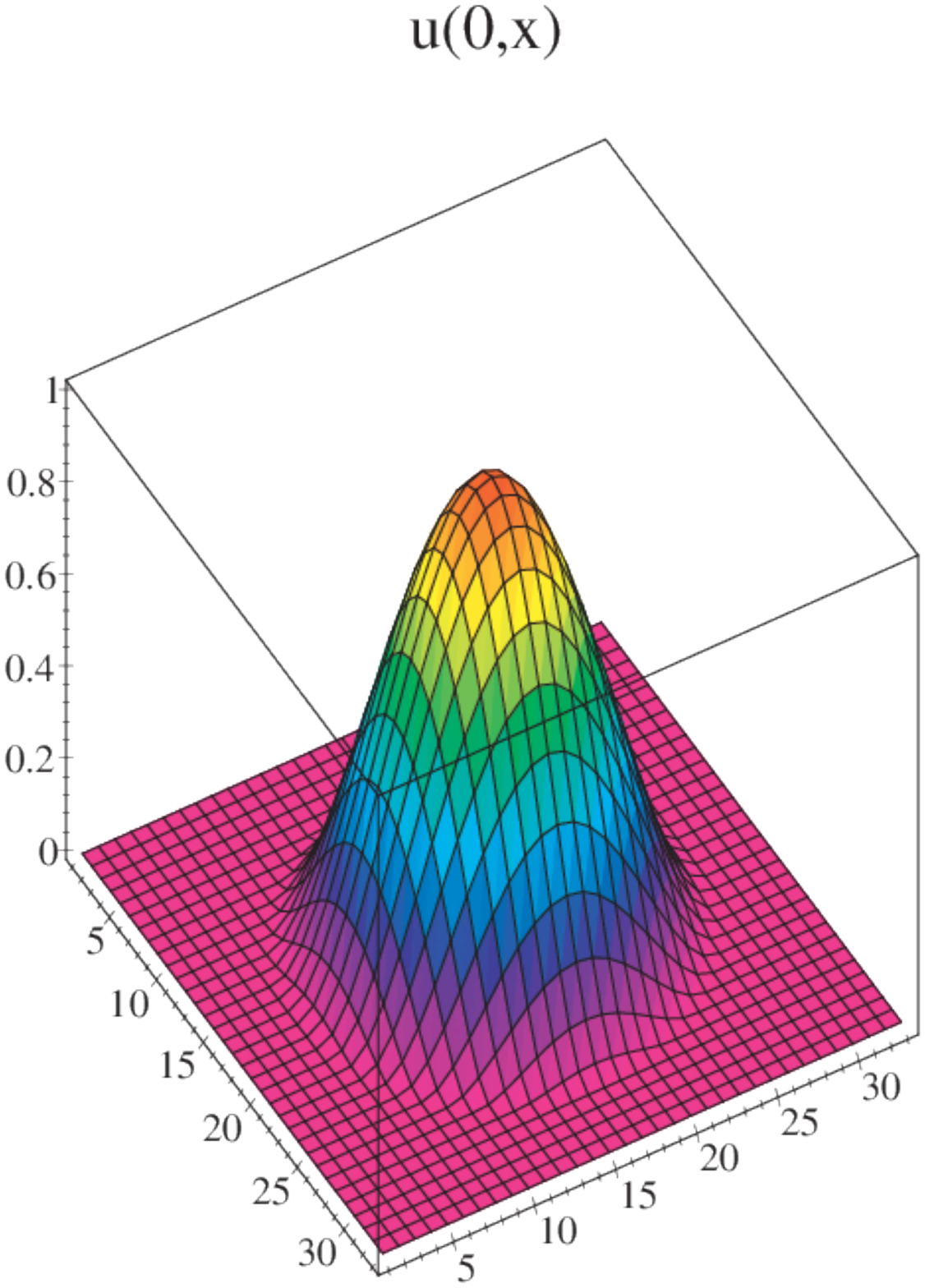} }
\vskip-1cm
{\bf Figure 3:} First choice of $f=u(T)$ and corresponding $u(0)$.
\vskip0.5cm

\centerline{ \epsfysize6cm \epsfxsize7cm \epsfbox{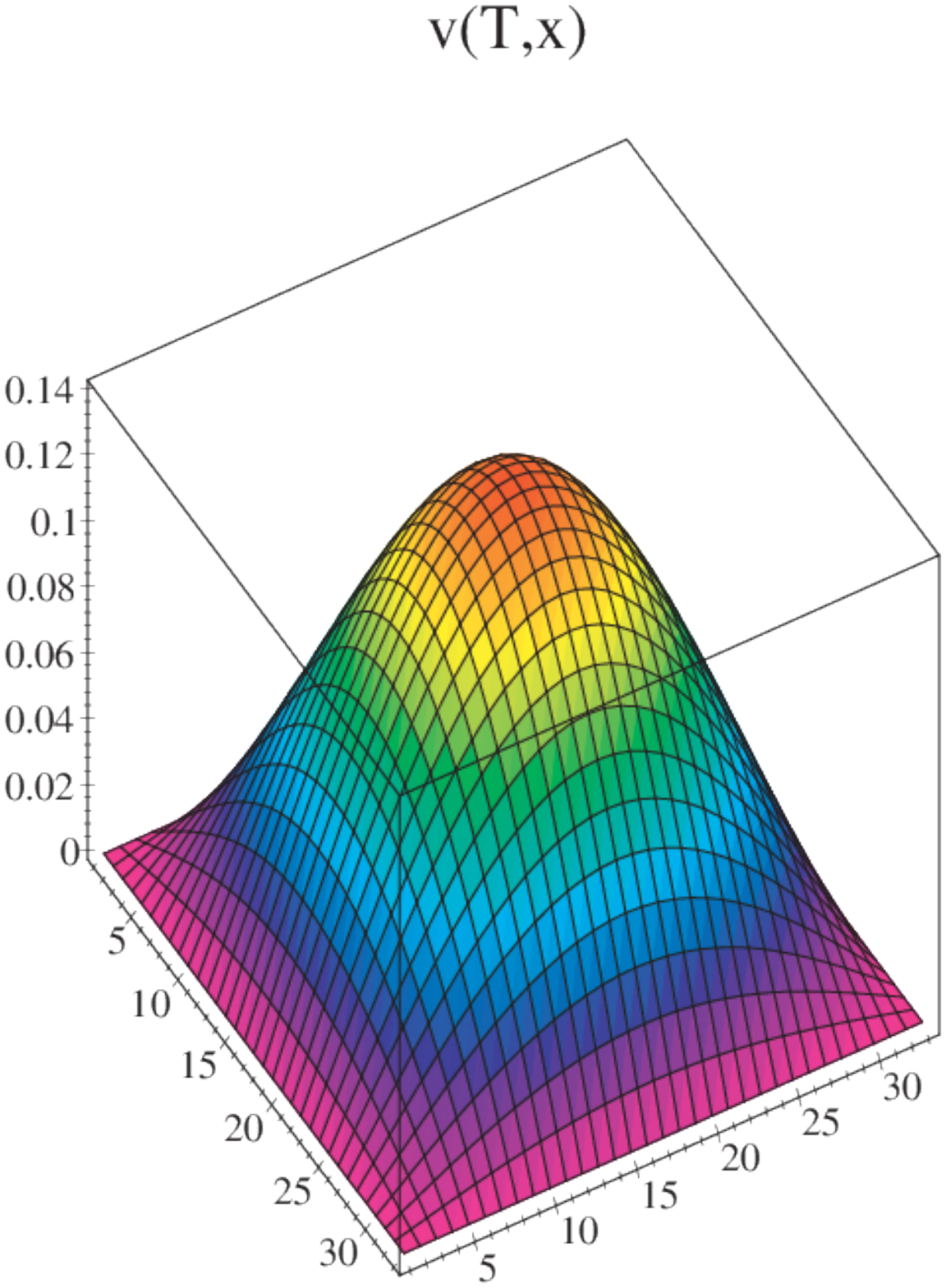}
             \epsfysize6cm \epsfxsize7cm \epsfbox{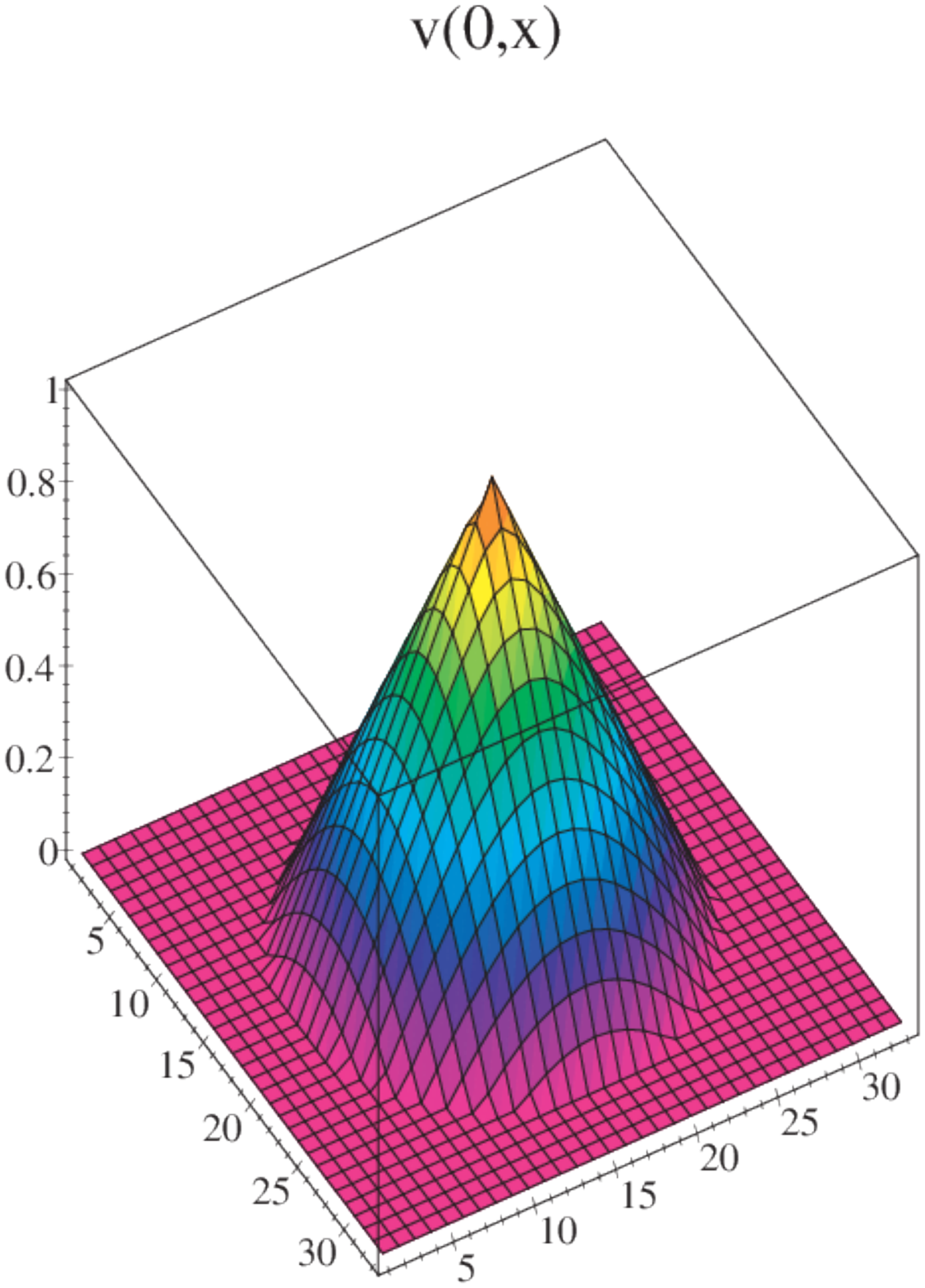} }
\vskip-1cm
{\bf Figure 4:} Second choice of $f=v(T)$ and corresponding $v(0)$.
\bigskip

\noindent  In Figure~5 the evolution of the iteration error for both problems 
is shown after $10$, $10^3$ and $10^5$ steps.
\newpage

\centerline{ \epsfysize6cm \epsfxsize7cm \epsfbox{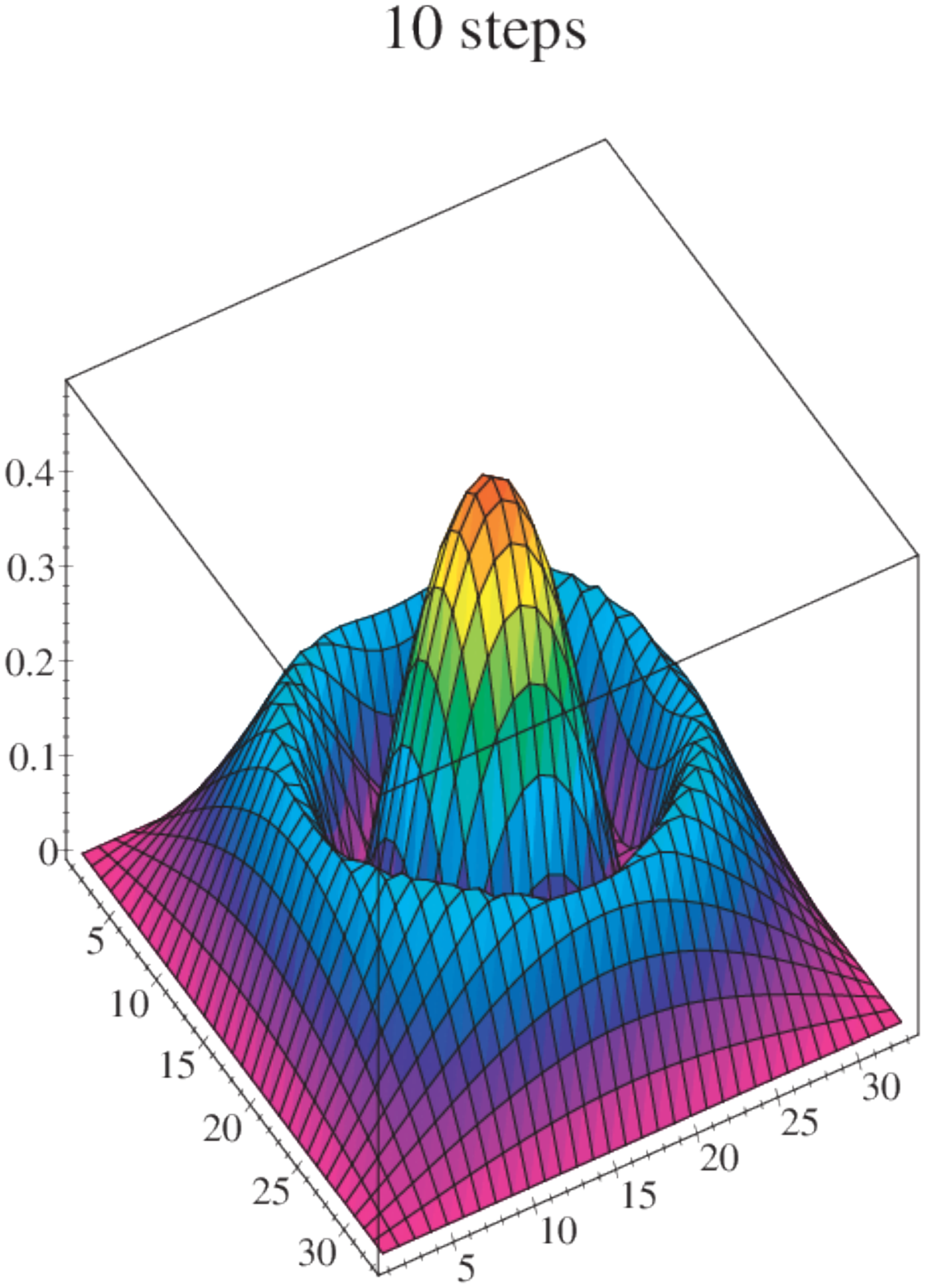}
             \epsfysize6cm \epsfxsize7cm \epsfbox{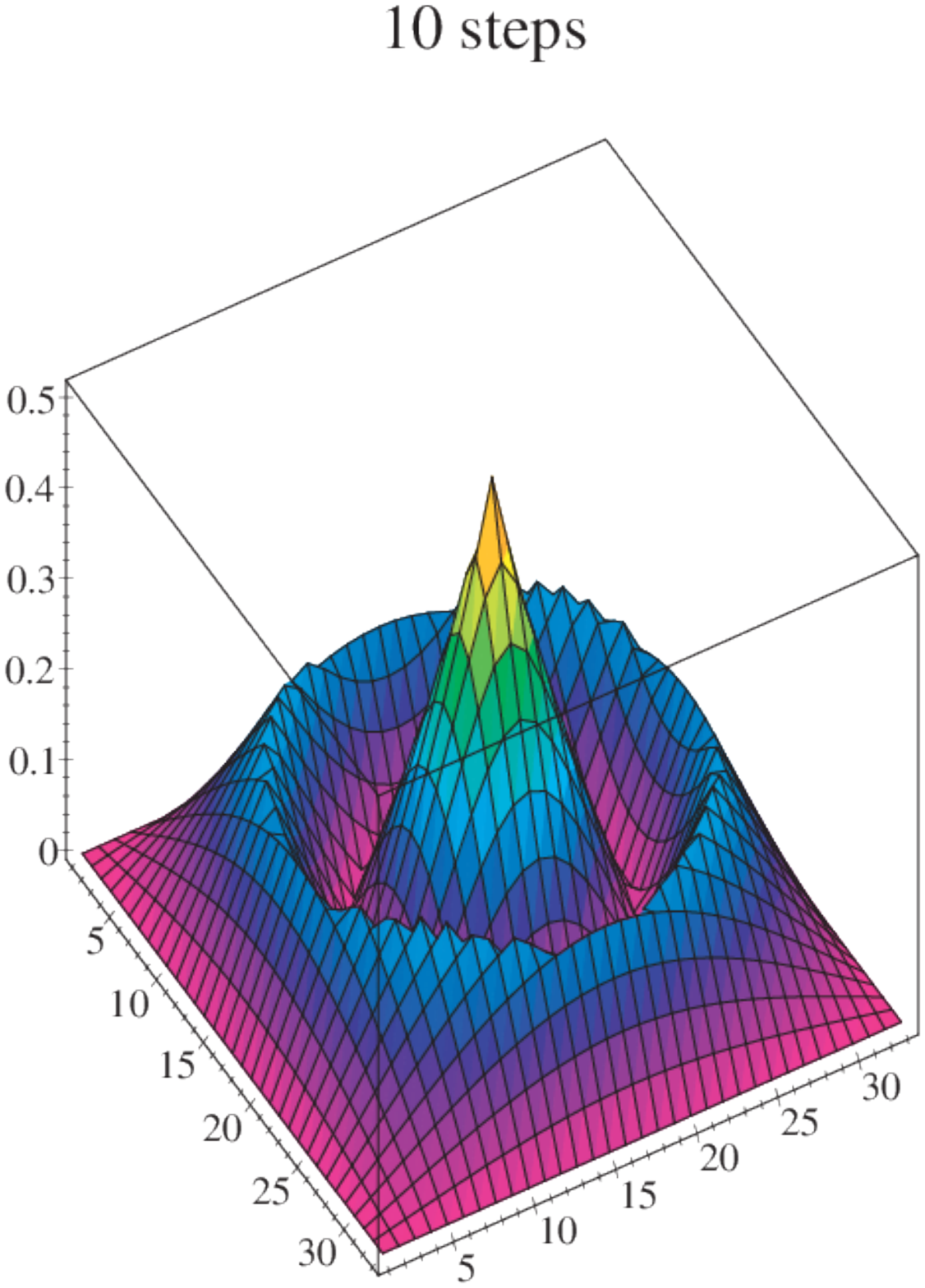} }

\centerline{ \epsfysize6cm \epsfxsize7cm \epsfbox{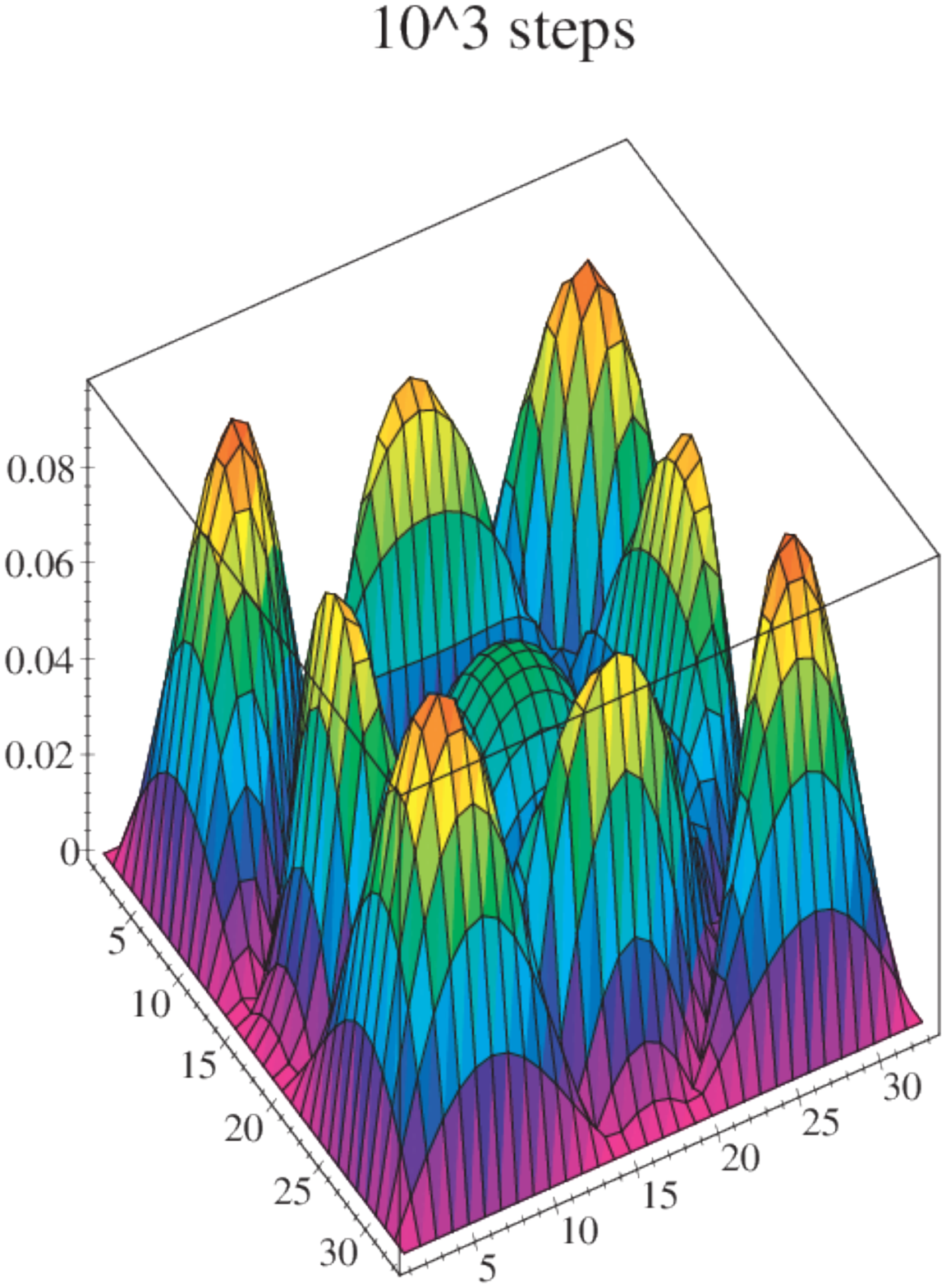}
             \epsfysize6cm \epsfxsize7cm \epsfbox{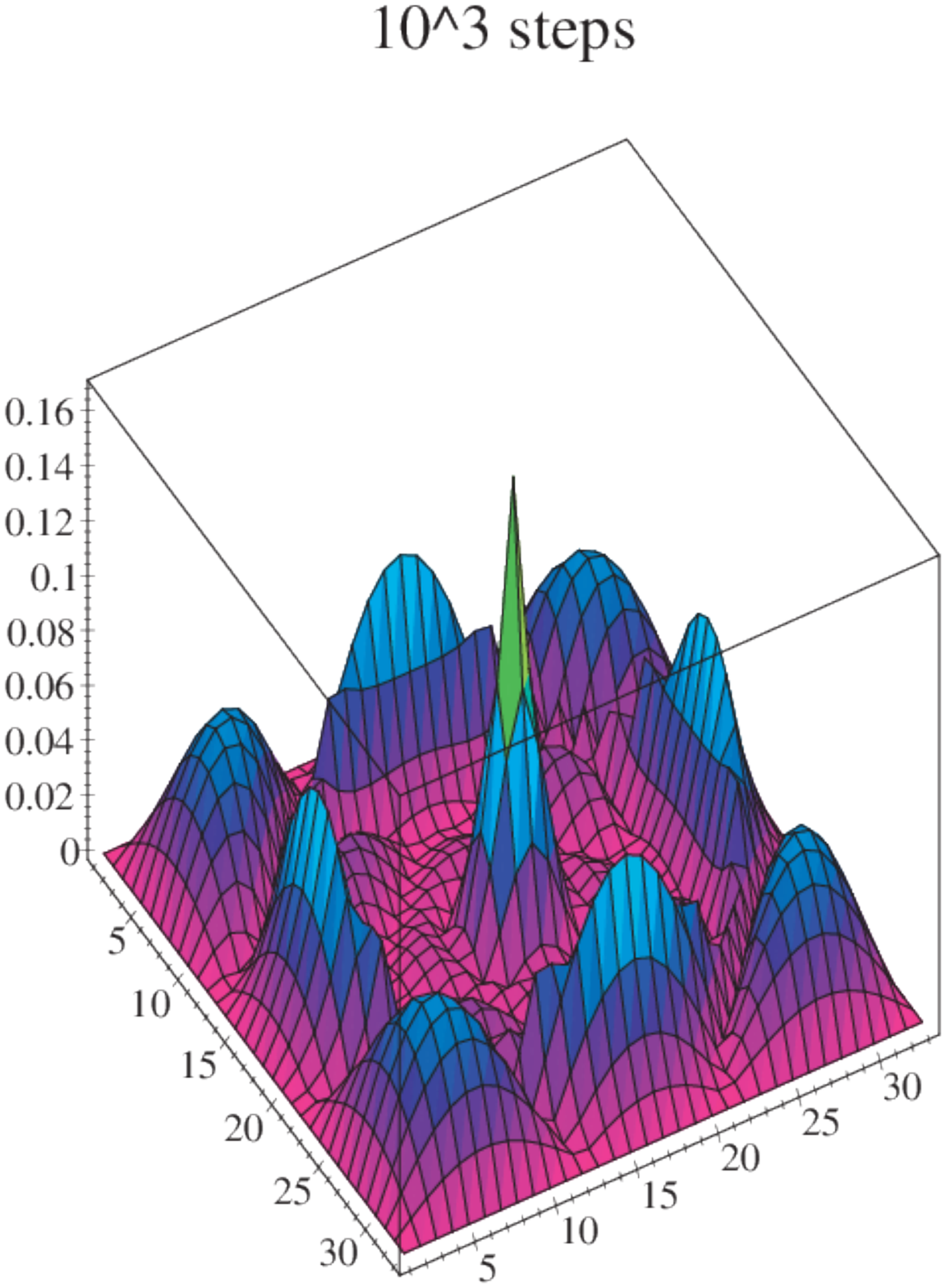} }

\centerline{ \epsfysize6cm \epsfxsize7cm \epsfbox{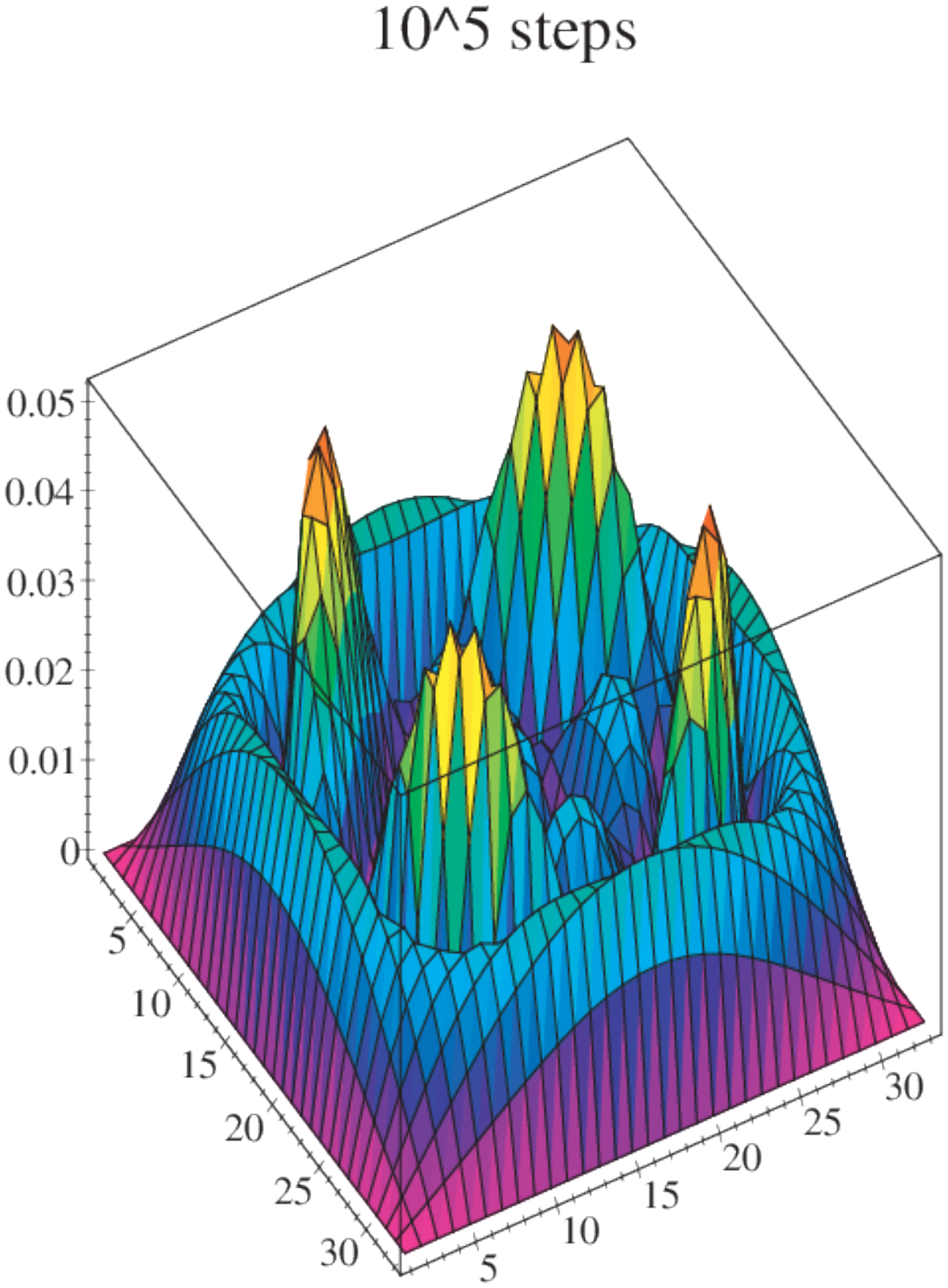}
             \epsfysize6cm \epsfxsize7cm \epsfbox{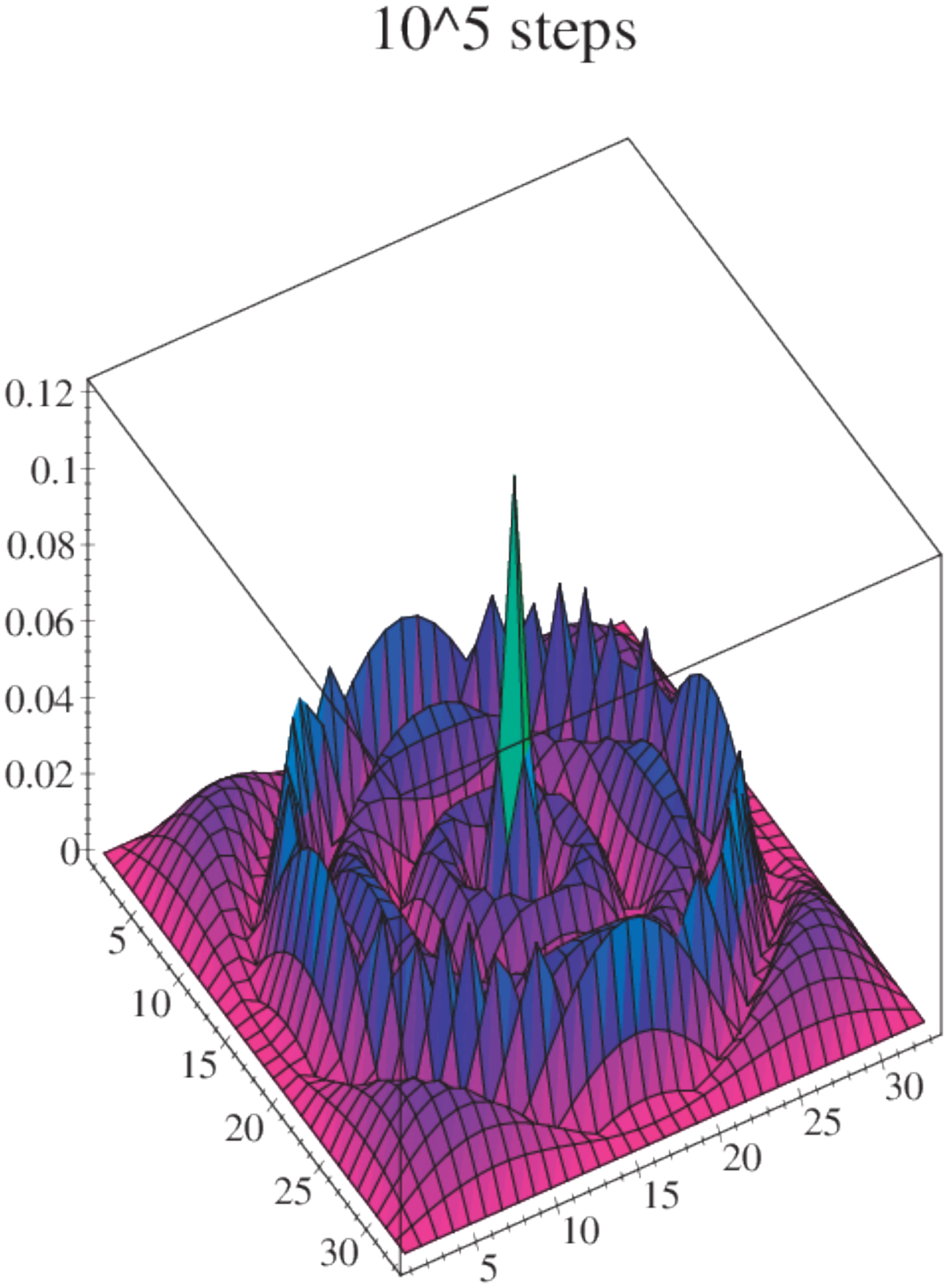} }
\vskip-1cm

{\bf Figure 5:} Evolution of the error for the boundary conditions $u(T)$ 
and $v(T)$ respectively.

\noindent  In Table~2 we present the evolution of the relative error in the 
$L^2$--norm for both boundary conditions $u(T)$ and $v(T)$.

\begin{center} \begin{tabular}{@{}cccccc} \hline\hline
{} & $10$ steps & $10^2$ steps & $10^3$ steps & $10^4$ steps & $10^5$ steps \\
\hline\hline
$f = u(T)$ & 53.5\%  & 37.9\% & 15.0\% & 12.1\% &  6.6\% \\ \hline
$f = v(T)$ & 46.1\%  & 32.2\% & 10.6\% &  9.1\% &  6.7\% \\ \hline\hline
\end{tabular} \\[1ex]
{\bf Table 2:} Evolution of the relative error in the $L^2$--norm.
\end{center}

\end{examp}

\section*{\large \bf References}
\begin{description}
\item[{[1]}] {\sc G. Bastay,} {\it Iterative Methods for Ill-Posed
     Boundary value Problems,} Lin\-k\"oping Studies in Science and Technology,
     Dissertations n. {\bf 392}, Link\"o\-ping, 1995

\item[{[2]}] {\sc M. Jourhmane and A. Nachaoui,} {\it A Relaxation
     Algorithm for Solving a Cauchy--Problem,} Preliminary Proceedings--Vol
     {\bf 2}, 2nd Intern. Confer. on Inverse Problems in Engineering: Theory
     and Practice, Le Croisic, 1996

\item[{[3]}] {\sc V.A. Kozlov, V.G. Maz'ya and A.V. Fomin,} {\it An
     iterative method for solving the Cauchy problem for elliptic equations,}
     Comput. Math. Phys., {\bf 31} (1991), no. {\bf 1}, 45 -- 52

\item[{[4]}] {\sc V.A. Kozlov and V.G. Maz'ya,} {\it On iterative
     procedures for solving ill-posed boundary value problems that preserve
     differential equations,} Lenin\-grad Math. J., {\bf 1} (1990), no.
     {\bf 5}, 1207 -- 1228

\item[{[5]}] {\sc A. Leit\~ao,} {\it Ein Iterationsverfahren f\"ur
     elliptische Cauchy--Probleme und die Ver\-kn\"upfung mit der
     Backus--Gilbert Methode,} Dissertation, FB Mathematik, J.W.
     Goethe--Universit\"at, Frankfurt am Main, 1996

\item[{[6]}] {\sc A. Leit\~ao,} {\it An Iterative Method for Solving
     Elliptic Cauchy Problems,} Numer. Funct. Anal. Optimization, {\bf 21}
     (2000), no. {\bf 5--6}, 715 -- 742

\item[{[7]}] {\sc G.M. Vainikko,} {\it Regularisierung nichkorrekter
      Aufgaben,} Preprint n. 200, Universit\"at Kaiserslautern, 1991
\end{description}
\vspace{.4in}
\begin{flushleft}
 Department of Mathematics  \\
 Federal University of Santa Catarina \\
 P.O. Box 476 \\
 88010-970, Florian\'opolis, BRAZIL \\
{ \em e-mail: aleitao@mtm.ufsc.br}
\end{flushleft}

\end{document}